# Об эффективных численных методах решения задач энтропийно-линейного программирования

*А.В. Гасников[1,2], Е.В. Гасникова[2], Ю.Е. Нестеров[3], А.В. Чернов[2]*

[1] Институт проблем передачи информации им. А.А. Харкевича Российской академии наук.
127051, Россия, г. Москва, Большой Каретный переулок, д.19 стр. 1
[2] Кафедра математических основ управления Факультета управления и прикладной математики Национального исследовательского Университета «Московский физико-технический институт» (НИУ МФТИ).
141700, Россия, Московская область, г. Долгопрудный, Институтский переулок, д. 9
[3] Center for Operation Research and Econometrics Université Catholique de Louvain.
Voie du Roman Pays 34, L1.03.01 - B-1348 Louvain-la-Neuve (Belgium)

**Аннотация**

В различных приложениях возникают задачи энтропийно-линейного программирования (ЭЛП). Эти задачи обычно записываются как задачи максимизации энтропии (минимизации минус энтропии) при аффинных ограничениях. В работе приводятся новые численные методы решения задач ЭЛП. Устанавливаются точные оценки скоростей сходимости предложенных методов. Изложенный в статье подход применим к более широкому классу задач минимизации сильно выпуклых функционалов при аффинных ограничениях.

**Ключевые слова:** энтропийно-линейное программирование, быстрый градиентный метод, регуляризация, двойственная задача.

## 1. Введение

Данная статья посвящена разработке эффективных численных методов решения задач энтропийно-линейного программирования (ЭЛП) [1] и получению оценок их скоростей сходимости. Имеется огромное число приложений, в которых возникают задачи ЭЛП [2 – 9].

Основная идея предлагаемого в статье подхода заключается в решении специальным образом регуляризованной (по Тихонову) двойственной задачи к задаче ЭЛП. Зная решение двойственной задачи, по явным формулам можно найти решение исходной задачи. В статье приводятся точные формулы, определяющие сколько достаточно сделать итераций быстрого градиентного метода для регуляризованной двойственной задачи, чтобы гарантировано (с заданной точностью) восстановить решение исходной задачи. Насколько нам известно, для задач ЭЛП такого рода формулы выписываются впервые. Собственно, по известным нам ос-тальным методам решения задач ЭЛП в основном имеются только теоремы о сходимости [1, 2, 10] (без оценок скорости сходимости).

Отметим также, что полученные в данной статье оценки при определенных условиях лучше имеющихся нижних оценок. Поясним сказанное.

Известно, что задача поиска такого $x^* \in \mathbb{R}^n$, что $Ax^* = b$ сводится к задаче выпуклой гладкой оптимизации

$$f(x) = \|Ax - b\|_2^2 \to \min_x.$$

Нижние оценки (при $k \le n$) для такой задачи приводят к оценке:

$$f(x_k) \ge \Omega\left(L_x R_x^2 / k^2\right),$$

т.е. существует такое $\chi$, что при $1 \le k \le n$ имеет место неравенство $f(x_k) \ge \chi L_x R_x^2 / k^2$. Здесь $k$ – число умножений матрицы $A$ на столбец/строку, $x_k$ – то, что можно получить на основе результатов этих умножений для самого лучшего метода на классе всевозможных задач такого типа; $L_x = \sigma_{\max}(A) = \lambda_{\max}(A^T A)$, $R_x = \|x^*\|_2 = \|A^T(AA^T)^{-1}b\|_2$. Откуда следует, что только при $k \ge \Omega\left(\sqrt{L_x} R_x / \varepsilon\right)$ можно гарантировать выполнение неравенства $f(x_k) \le \varepsilon^2$, т.е. $\|Ax_k - b\|_2 \le \varepsilon$. Детали см. в [стр. 264–274, 11]. Далее в статье будет показано (см., например, теорему 2), что при достаточно общих условиях можно улучшить приведенную оценку, попутно добившись, чтобы полученный $x_k$ был бы приближенным решением задачи ЭЛП, т.е. не только удовлетворял бы "возмущенным" аффинным ограничениям $\|Ax_k - b\|_2 \le \varepsilon$ в задаче ЭЛП.

Структура статьи следующая. В п. 2 приведен простейший пример возникновения задачи ЭЛП. В п. 3 рассматривается постановка задачи ЭЛП, приводятся основные определения. В п. 4 приводятся необходимые вспомогательные результаты. В п. 5 формулируются основные результаты работы. В п. 6 эти результаты обсуждаются. В п. 7 делаются заключительные замечания.

## 2. Парадокс Эренфестов

Следуя учебнику [12] опишем известный в физике парадокс Эренфестов [13] в немного вольной трактовке.

Рядом стоят две собаки с номерами *1* и *2*. На собаках как-то расположились $M \gg 1$ блох. Каждая блоха в промежутке времени $[t, t+h)$ с вероятностью $\lambda h + o(h)$ ($\lambda = 1$) независимо от остальных перескакивает на соседнюю собаку. Пусть в начальный момент все блохи собрались на собаке с номером *1*. Тогда для всех $t \geq \chi M$ ($\chi \sim 10$)

$$P\left(\frac{|n_1(t) - n_2(t)|}{M} \leq \frac{5}{\sqrt{M}}\right) \geq 0.99,$$

где $n_1(t)$ – число блох на первой собаке в момент времени $t$, а $n_2(t)$ – на второй (случайные величины). То есть относительная разность числа блох на собаках будет иметь порядок малости $O(1/\sqrt{M})$ на больших временах ($t \geq \chi M$). Далее мы поясним связь этого результата с принципом максимума энтропии, который мы будем записывать в виде принципа минимума минус энтропии.

Описанная выше марковская динамика имеет закон сохранения числа блох: $n_1(t) + n_2(t) \equiv M$, и это будет единственный закон сохранения. Стационарная мера имеет вид (теорема Санова):

$$\nu(n_1, n_2) = \nu(c_1 M, c_2 M) = M! \frac{(1/2)^{n_1}}{n_1!} \frac{(1/2)^{n_2}}{n_2!} = C_M^{n_1} 2^{-M} \simeq \frac{2^{-M}}{\sqrt{2\pi c_1 c_2}} \exp(-M \cdot H(c_1, c_2)),$$

где $H(c_1, c_2) = \sum_{i=1}^{2} c_i \ln c_i$, $c_i$ – концентрация блох на собаки $i$ в равновесии (т.е. при $t \to \infty$). Кстати, сказать, из такого вида стационарной меры, можно получить, что если в начальный момент все блохи находились на одной собаке, то математическое ожидание времени первого возвращения макросистемы в такое состояние будет порядка $2^M$. Равновесие данной макросистемы естественно определять как состояние, в малой окрестности которого сконцентрирована стационарная мера (принцип максимума энтропии Больцмана–Джейнса)

$$c^* = \begin{pmatrix} 1/2 \\ 1/2 \end{pmatrix} = \arg\min_{\substack{c_1 + c_2 = 1 \\ c \geq 0}} H(c).$$

Поиск равновесия привел к необходимости решения задачи ЭЛП. Этот же результат можно было получить и при другом порядке предельных переходов (обратном к рассмотренному выше порядку: $t \to \infty$, $M \to \infty$). А именно, сначала считаем, что при $t = 0$ существует предел $c_i(t) \stackrel{\text{п.н.}}{=} \lim_{M \to \infty} n_i(t)/M$. Тогда (теорема Т. Куртца) этот предел существует при любом $t > 0$, причем $c_1(t)$, $c_2(t)$ – детерминированные (не случайные) функции, удовлетворяющие СОДУ

$$\frac{dc_1}{dt} = \lambda \cdot (c_2 - c_1),$$

$$\frac{dc_2}{dt} = \lambda \cdot (c_1 - c_2).$$

Глобально устойчивым положением равновесия этой СОДУ будет $c^*$, а $H(c)$ – функция Ляпунова этой СОДУ (убывает на траекториях СОДУ, и имеет минимум в точке $c^*$). Все это можно перенести на общие модели стохастической химической кинетики с (обобщенным) условием детального баланса [8, 14 – 16].

Приведенный парадокс Эренфестов является, пожалуй, простейшим примером того как возникают задачи ЭЛП при поиске равновесий макросистем. Подробнее об этом можно прочитать, например, в [8]. В частности, в большинстве реальных приложений возникают задачи ЭЛП, в которых не "две собаки" и не "один закон сохранения числа блох", а на много больше (см., например, [17]), и возникающая задача ЭЛП не решается по явным формулам. Требуется разработка эффективных численных методов.

### 3. Задача энтропийно-линейного программирования

Рассматривается задача ЭЛП [1]

$$f(x) = \sum_{i=1}^{n} x_i \ln(x_i/\xi_i) \to \min_{x \in S_n(1); Ax=b}, \qquad (1)$$

где $S_n(1) = \left\{ x \in \mathbb{R}^n : x_i \geq 0,\ i=1,\ldots,n,\ \sum_{i=1}^{n} x_i = 1 \right\}$ – единичный симплекс [1] в $\mathbb{R}^n$, $\xi \in \operatorname{ri} S_n(\Xi)$.

Мы считаем матрицу $A$ разреженной с не более чем $s \ll n$ элементами в каждой строке. Число строк в матрице $A$ есть $m \ll n$. Построим двойственную задачу

$$\min_{x \in S_n(1);\, Ax=b} \sum_{i=1}^{n} x_i \ln(x_i/\xi_i) = \min_{x \in S_n(1)} \max_{\lambda \in \mathbb{R}^m} \left\{ \sum_{i=1}^{n} x_i \ln(x_i/\xi_i) + \langle \lambda, b - Ax \rangle \right\} =$$

$$= \max_{\lambda \in \mathbb{R}^m} \min_{x \in S_n(1)} \left\{ \sum_{i=1}^{n} x_i \ln(x_i/\xi_i) + \langle \lambda, b - Ax \rangle \right\} = \max_{\lambda \in \mathbb{R}^m} \left\{ \langle \lambda, b \rangle - \ln\left( \sum_{i=1}^{n} \xi_i \exp\left( \left[ A^T \lambda \right]_i \right) \right) \right\}.$$

Таким образом, двойственная задача имеет вид

$$\varphi(\lambda) = \langle \lambda, b \rangle - \ln\left( \sum_{i=1}^{n} \xi_i \exp\left( \left[ A^T \lambda \right]_i \right) \right) \to \max_{\lambda \in \mathbb{R}^m}. \qquad (2)$$

Решения прямой и двойственной задачи связаны:

---

[1] Более общая ситуация, когда $x \in S_n(\Lambda)$, может быть сведенная к $x \in S_n(1)$ с помощью замены переменных. При этом $\tilde{x} := x/\Lambda$, $\tilde{\xi} := \xi/\Lambda$, $\tilde{b} := b/\Lambda$, $\tilde{A} := A$, $\tilde{f} := f\Lambda$.

$$x_i(\lambda) = \frac{\xi_i \exp\left(\left[A^T\lambda\right]_i\right)}{\sum_{k=1}^{n} \xi_k \exp\left(\left[A^T\lambda\right]_k\right)}.$$

**Определение.** *Под $(\varepsilon_f, \varepsilon)$-решением задачи (1) будем понимать такой вектор $x$, что*

$$f(x) - f^* \leq \varepsilon_f, \quad \|Ax - b\|_2 \leq \varepsilon,$$

*где* $f^* = \min\limits_{x \in S_n(1); Ax = b} f(x)$.

**Лемма 1.** *Пусть* $-\langle \lambda, \nabla\varphi(\lambda)\rangle \leq \varepsilon_f$, $\|\nabla\varphi(\lambda)\|_2 \leq \varepsilon$, *тогда $x(\lambda)$ – $(\varepsilon_f, \varepsilon)$-решение задачи (1).*

**Доказательство.** Соотношение $\nabla\varphi(\lambda) = b - Ax(\lambda)$ следует из представления двойственной задачи $\varphi(\lambda) = \min\limits_{x \in S_n(1)} \{f(x) + \langle \lambda, b - Ax\rangle\}$ и теоремы Демьянова–Рубинова–Данскина.

Обозначим через $x^*$ решение задач (1). Тогда

$$f(x(\lambda)) + \langle \lambda, b - Ax(\lambda)\rangle \leq f(x^*) + \langle \lambda, b - Ax^*\rangle = f(x^*) = f^*.$$

Откуда

$$f(x(\lambda)) = f(x^*) - \langle \lambda, b - Ax(\lambda)\rangle = f^* - \langle \lambda, \nabla\varphi(\lambda)\rangle \leq f^* + \varepsilon_f. \bullet$$

**Лемма 2 (см. [18]).** *Имеет место следующее неравенство*

$$\|\nabla\varphi(\lambda_2) - \nabla\varphi(\lambda_1)\|_2 \leq L\|\lambda_2 - \lambda_1\|_2,$$

*где* $L = \max\limits_{1 \leq i \leq n} \left\|[A]^{(i)}\right\|_2^2$, $[A]^{(i)}$ – *$i$-й столбец матрицы $A$.*

К сожалению, несмотря на установленную в лемме 2 гладкость функционала (2), мы имеем лишь довольно грубые оценки его сильной вогнутости. Быстрый градиентный метод (БГМ) [19] (стартующий с $\lambda = 0$) для задачи (2) после $\mathrm{O}\left(\sqrt{LR^2/\tilde{\varepsilon}}\right)$ итераций, где $R = \|\lambda^*\|_2$ – размер решения задачи (2), $\tilde{\varepsilon} = \min\left\{\varepsilon^2/(2L), \varepsilon_f^2/(2LR^2)\right\}$ – точность решения задачи (2), га-рантирует в виду леммы 1 и следующих неравенств ($\varphi^*$ – максимальное значение функционала в двойственной задаче)

$$\frac{1}{2L}\|\nabla\varphi(\lambda)\|_2^2 \le \varphi^* - \varphi(\lambda), \quad \frac{|-\langle\lambda,\nabla\varphi(\lambda)\rangle|^2}{2LR^2} \le \frac{1}{2L}\|\nabla\varphi(\lambda)\|_2^2 \le \varphi^* - \varphi(\lambda),$$

что $x(\lambda) - (\varepsilon_f,\varepsilon)$-решение задачи (1). Итоговая оценка на число итераций примет вид

$$\mathrm{O}\left(\max\left\{LR/\varepsilon, LR^2/\varepsilon_f\right\}\right). \tag{3}$$

Эту оценку можно улучшить. Этому и будет посвящена последующая часть данной работы.

## 4. Вспомогательные результаты для регуляризованной двойственной задачи к задаче энтропийно-линейного программирования

Регуляризуем функционал (2) по А.Н. Тихонову:

$$\varphi_\delta(\lambda) = \varphi(\lambda) - \frac{\delta}{2}\|\lambda\|_2^2,$$

и вместо задачи (2) будем решать регуляризованную задачу

$$\varphi_\delta(\lambda) \to \max_{\lambda\in\mathbb{R}^m}.$$

Параметр $\delta$ будет оптимально подобран позже.

Следующие три леммы решают задачу восстановления $(\varepsilon_f,\varepsilon)$-решением задачи (1) по решению этой регуляризованной задачи. В отличие от нерегуляризованного случая здесь возникают некоторые технические места, требующие аккуратной проработки.

**Лемма 3.** *Имеют место следующие неравенства*

$$\|\nabla\varphi(\lambda)\|_2 \le \|\nabla\varphi_\delta(\lambda)\|_2 + \delta\|\lambda\|_2,$$

$$-\langle\lambda,\nabla\varphi(\lambda)\rangle \le \frac{L_\delta}{\delta}\left(\varphi_\delta^* - \varphi_\delta(\lambda)\right), L_\delta = L + \delta,$$

$$\frac{\delta}{2}\left\|\lambda_{\delta}^{*}\right\|_{2}^{2} \le \ln\left(\Xi\Big/\min_{i=1,\ldots,n}\xi_{i}\right) \stackrel{def}{=} \Delta_{\varphi},$$

где $\lambda_{\delta}^{*}$ – решение регуляризованной двойственной задачи, а $\varphi_{\delta}^{*}$ – значение функционала в регуляризованной двойственной задаче на этом решении.

**Доказательство.** Не очень тривиально лишь второе и третье неравенство. Второе неравенство следует из следующей выкладки

$$\varphi_{\delta}^{*} - \varphi_{\delta}(\lambda) \ge \frac{\left\|\nabla\varphi_{\delta}(\lambda)\right\|_{2}^{2}}{2L_{\delta}} = \frac{\left\|\nabla\varphi(\lambda) - \delta\lambda\right\|_{2}^{2}}{2L_{\delta}} \ge -\frac{\delta\langle\nabla\varphi(\lambda),\lambda\rangle}{L_{\delta}}.$$

Третье неравенство следует из следующей выкладки

$$\frac{\delta}{2}\left\|\lambda_{\delta}^{*}\right\|_{2}^{2} \le \varphi_{\delta}^{*} - \varphi_{\delta}(0) \le \varphi^{*} - \varphi(0) =$$

$$= \min_{x\in S_{n}(1);\, Ax=b}\sum_{i=1}^{n}x_{i}\ln(x_{i}/\xi_{i}) - \min_{x\in S_{n}(1)}\sum_{i=1}^{n}x_{i}\ln(x_{i}/\xi_{i}) \le \ln\left(\Xi\Big/\min_{i=1,\ldots,n}\xi_{i}\right).$$

Первое неравенство в этой цепочке следует из того, что $\nabla\varphi_{\delta}(\lambda_{\delta}^{*})=0$ и $\varphi_{\delta}(\lambda)$ – сильно вогнутая функция, с константой сильной вогнутости $\ge \delta$. ●

Для решения регуляризованной задачи воспользуемся БГМ для сильно вогнутых задач [19] ($\lambda_{0}=u_{0}=0$):

$$\begin{cases}\lambda_{k+1}=u_{k}+\dfrac{1}{L_{\delta}}\big(b-Ax(\lambda_{k})-\delta\lambda_{k}\big),\\[2mm] u_{k+1}=\lambda_{k}+\dfrac{\sqrt{L_{\delta}}-\sqrt{\delta}}{\sqrt{L_{\delta}}+\sqrt{\delta}}(\lambda_{k+1}-\lambda_{k}).\end{cases} \qquad(4)$$

**Лемма 4 (см. [19]).** *Имеют место следующие неравенства ($k=0,1,\ldots$)*

$$\varphi_{\delta}^{*}-\varphi_{\delta}(\lambda_{k}) \le 2\Delta_{\varphi}\exp\left(-k\sqrt{\frac{\delta}{L_{\delta}}}\right),$$

$$\left\|\nabla\varphi_{\delta}(\lambda_{k})\right\|_{2}^{2} \le 4L_{\delta}\Delta_{\varphi}\exp\left(-k\sqrt{\frac{\delta}{L_{\delta}}}\right),$$

$$\left\|\lambda_k - \lambda_\delta^*\right\|_2^2 \le \min\left\{\left\|\lambda_\delta^*\right\|_2^2, \frac{4\Delta_\varphi}{\delta}\exp\left(-k\sqrt{\frac{\delta}{L_\delta}}\right)\right\}.$$

Из лемм 3, 4 следует

**Лемма 5.** *Имеет место следующее неравенство*

$$\left\|\lambda_k\right\|_2 \le \left\|\lambda_k - \lambda_\delta^*\right\|_2 + \left\|\lambda_\delta^*\right\|_2 \le 2\left\|\lambda_\delta^*\right\|_2 \le \sqrt{8\Delta_\varphi/\delta}.$$

## 5. Основные результаты

Следствием лемм 1, 3 – 5 является следующая

**Теорема 1.** *Если выбрать* $\sqrt{\delta} \simeq \varepsilon/\sqrt{9\Delta_\varphi}$, *то после*

$$N \simeq \frac{\sqrt{9L\Delta_\varphi}}{\varepsilon}\ln\left(\frac{9L\Delta_\varphi}{\varepsilon_f \varepsilon^2}\right) \tag{5}$$

*итераций метода (4) получим такой* $\lambda_N$, *что* $x(\lambda_N) - (\varepsilon_f, \varepsilon)$-*решением задачи (1)*.

**Доказательство.** Из лемм 1, 3 – 5 имеем

$$\exp\left(-N\sqrt{\frac{\delta}{L_\delta}}\right) = \min\left\{\frac{\delta\varepsilon_f}{2L_\delta\Delta_\varphi}, \frac{\left(\varepsilon - \sqrt{8\Delta_\varphi\delta}\right)_+^2}{4L_\delta\Delta_\varphi}\right\}.$$

Откуда $\sqrt{2\varepsilon_f}\sqrt{\delta} = \varepsilon - \sqrt{8\Delta_\varphi}\sqrt{\delta}$, следовательно $\sqrt{\delta} = \varepsilon/\left(\sqrt{8\Delta_\varphi} + \sqrt{2\varepsilon_f}\right)$. Для упрощения формул будем выбирать $\delta$ немного не оптимально: $\sqrt{\delta} \simeq \varepsilon/\sqrt{9\Delta_\varphi}$. Тогда

$$N \simeq \frac{\sqrt{9L\Delta_\varphi}}{\varepsilon}\ln\left(\frac{18L\Delta_\varphi}{\varepsilon_f \varepsilon^2}\right). \bullet$$

В отличие от полученной ранее оценки (3), в оценку (5) не входит потенциально большой размер $R = \left\|\lambda^*\right\|_2$. О том, что в приложениях этот размер, действительно, может быть большим, говорят результаты численных экспериментов [20]. Оценка (5) доминирует оценку

(3) даже для небольших значений $R$. Естественным образом возникает вопрос: можно ли улучшить оценку (5) если известно, что $R$ – не очень большое число? Разберем эту ситуацию.

Обозначим через $R_\delta = \max_{k=0,1,...} \|\lambda_k\|_2$. Лемма 5 дает оценку на $R_\delta$, но эта оценка может быть сильно завышена, поэтому завышенной может получиться и оценка (5). Представим, что нам известно значение $R$, которое мажорирует $R_\delta$ (равномерно по $\delta \geq 0$). Выберем оптимально $\delta$ исходя из известного значения $R$. Для этого в виду лемм 3, 5 нужно решить задачу (см. также доказательство теоремы 1 с $R = \sqrt{8\Delta_\varphi/\delta}$)

$$\max\left\{2\varepsilon_f \delta, (\varepsilon - R\delta)_+^2\right\} \to \min_{\delta \geq 0}.$$

Эту задачу можно решить явно

$$\sqrt{\delta} = \frac{\varepsilon}{\sqrt{\varepsilon_f/2} + \sqrt{\varepsilon_f/2 + R\varepsilon}}.$$

Однако нам будет удобнее взять

$$\sqrt{\delta} \simeq \frac{\varepsilon}{2\sqrt{\varepsilon_f/2 + R\varepsilon}}. \tag{6}$$

Такая замена не сильно скажется на оценках необходимого числа итераций, но несколько упростит формулы.

**Теорема 2.** *Если выбрать $\delta$ согласно формуле (6), то после*

$$N \simeq \sqrt{\frac{2L \cdot (\varepsilon_f + 2R\varepsilon)}{\varepsilon^2}} \ln\left(\frac{4L\Delta_\varphi \cdot (\varepsilon_f + 2R\varepsilon)}{\varepsilon_f \varepsilon^2}\right) \tag{7}$$

*итераций метода (4) получим такой $\lambda_N$, что $x(\lambda_N) - (\varepsilon_f, \varepsilon)$-решением задачи (1).*

*Основной вклад в оценку общей трудоемкости вносит расчет градиента в (4)*

$$\nabla \varphi_\delta(\lambda) = b - Ax(\lambda) - \delta\lambda,$$

*требующий* $\mathrm{O}(n+sm)$ *арифметических операций. Общие затраты будут*

$$\mathrm{O}\left((n+sm)\sqrt{\frac{2L\cdot(\varepsilon_f+2R\varepsilon)}{\varepsilon^2}}\ln\left(\frac{4L\Delta_\varphi\cdot(\varepsilon_f+2R\varepsilon)}{\varepsilon_f\varepsilon^2}\right)\right)$$

Доказательство теоремы 2 аналогично доказательству теоремы 1.

## 6. Обсуждение основных результатов

Оценка (7) доминирует оценку (5) для не очень больших значений $R$. Для очень больших значений $R$ оценки (5) и (7) неплохо соответствуют друг другу. Покажем это.

**Лемма 6 (см. [21]).** *Пусть в задаче выпуклой оптимизации* $f(x) \to \min\limits_{Ax=b,\,x\in Q}$ *функция* $f(x)$ *обладает ограниченной вариацией на множестве* $Q$: $\max\limits_{x,y\in Q}(f(x)-f(y))\le\Delta$. *Предположим, что* $B_2(0,r)$ – *евклидов шар (в двойственном пространстве* $\lambda$ – *множитель Лагранжа к ограничению* $Ax=b$) *радиуса* $r$ *с центром в точке 0 полностью принадлежит множеству* $\Xi_{b,A}=\{\lambda:\lambda=b-Ax,\,x\in Q\}$. *Тогда имеет место следующая оценка на размер решения двойственной задачи:* $\|\lambda^*\|_2\le\Delta/r$.

Для задачи (1) $\Delta=\Delta_\varphi$, $Q=S_n(1)$. Нужно оценить $r$, чтобы можно было воспользоваться этой леммой. Следующую идею (слейтеровской релаксации) сообщил нам А.С. Немировский. Можно довольно грубо оценить снизу размер вписанного шар $r$. Для этого "подменим" исходную задачу немного другой, в которой вектор $b$ заменяется на такой вектор $b_\varepsilon$, что $\|b_\varepsilon-b\|_2\le\varepsilon$ и $B_2(0,\varepsilon)\subset\Xi_{b_\varepsilon,A}$. Тогда, для задачи с вектором $b_\varepsilon$ мы имеем оценку $r\ge\varepsilon$. Но $(\varepsilon_f,\varepsilon)$-решение этой обновленной задачи гарантировано будет $(\varepsilon_f,2\varepsilon)$-решением исходной задачи. Таким образом, в самом плохом случае можно считать, что $r\simeq\varepsilon$. Используя это наблюдение несложно сопоставить теоремы 1, 2 в случае, когда $R$ очень большое и $\varepsilon_f=\mathrm{O}(R\varepsilon)$.

Остается только одна проблема: неизвестность $R$. Хотя в шаг БГМ (4) $R$ явно не входит, $R$ входит в $\delta$ (см. формулу (6)), от которого размер шага уже зависит. Также $R$ входит в критерий останова метода (7), заключающегося в выполнении точно рассчитанного числа итераций. Далее (следуя [21, 22]) будет описана процедура рестартов, позволяющая справиться с этой сложностью.

Полагаем $R = R_0 = 100$ (при $n \gg 10^4$, $m \gg 10^2$ в различных численных экспериментах [20] получалось $R \gg 10^3$) делаем предписанное этому $R$ число итераций (7) и проверяем критерий останова – лемма 1. Заметим, что в этом критерии останова не надо знать $\varphi^*$. Если этот критерий не выполняется, то полагаем $R := 4R$ и повторяем процедуру. Через не более чем $\lceil \log_4(4R/R_0) \rceil$ перезапусков мы остановимся.

Поясним, почему была выбрана именно константа 4. Оптимально выбирать такой коэффициент $\beta$ ($R := \beta R$), который доставляет минимум следующему выражению $\beta\sqrt{R/R_0}/(\sqrt{\beta}-1)$, получающемуся при оценке сверху общего количества сделанных итераций с учетом всех необходимых перезапусков (считаем $R/R_0 \gg 1$, $\varepsilon_f = \mathrm{O}(R\varepsilon)$).[2] Решением этой задачи будет $\beta = 4$.

## 7. Заключительные замечания

Приведенные в статье результаты могут быть перенесены на другие сильно выпуклые функционалы. Выбор в качестве функционала минус энтропии лишь один из возможных вариантов. Напомним, что согласно неравенству Пинскера минус энтропия – сильно выпуклая функция в 1-норме с константой сильной выпуклости 1. При этом не так важно, чтобы решение прямой и двойственной задачи были связаны явными формулами, как это имеет место в данной статье, и, в целом, довольно типично для сепарабельных функционалов. Достаточно, чтобы имела место сильная выпуклость функционала (см. пример 4 [23]) или его сепарабельность вместе с ограничениями. Это, в частности, обеспечивает возможность с геометрической скоростью сходимости находить приближенную зависимость $x(\lambda)$.

---

[2] Действительно, при $\varepsilon_f = \mathrm{O}(R\varepsilon)$ теоремы 1 и 2 дают одинаковую зависимость $N \sim \sqrt{R}$. Таким образом, общее число сделанных итераций при $k$ рестартах пропорционально $N \sim \sqrt{\beta^{k+1}}/(\sqrt{\beta}-1)$, где натуральное число $k$ (большое при $R/R_0 \gg 1$) определяется из условия $\beta^{k-1}R_0 \leq R \leq \beta^k R_0$, т.е. $\sqrt{\beta^{k-1}} \leq \sqrt{R/R_0}$. Используя это соотношение, можно оценить сверху $\sqrt{\beta^{k+1}}/(\sqrt{\beta}-1)$ как $\beta\sqrt{R/R_0}/(\sqrt{\beta}-1)$.

Если размеры матрицы в аффинных ограничениях настолько большие, что умножение такой матрицы на столбец/строку не выполнимо за разумное время (напомним, что современный процессор может делать до миллиарда операций с плавающий точкой в секунду), то требуется другая техника решения отмеченных задач. В качестве одного из вариантов альтернативного подхода, укажем на идею седлового представления задачи (как правило, это требует и некоторого пересмотра самой постановки задачи), и последующее решение седловой задачи рандомизированным вариантом метода зеркального спуска или проксимальным зеркальным методом А.С. Немировского [24]. При этом рандомизация возникает на этапе расчета стохастического градиента, т.е. на этапе умножения матрицы на столбец/строку. Выгода от такой рандомизации – существенное удешевление стоимости одной итерации, а плата за это – увеличение числа итераций. К сожалению, такая рандомизация приводит к заметным результатам лишь тогда, когда в прямом и двойственном пространствах переменные "живут" на 1-шарах (или симплексах), что имеет место в лишь в очень небольшом числе приложений (типа задач поиска селектора Данцига [24]). В других случаях рандомизация также используется, но эффект, как правило, заметно меньше (см. концовку работы [24], п. 6.5 обзора [25], новую работу [26]).

Интересно было бы сравнить предложенные в статье методы с другими известными численными методами решения задач ЭЛП [1, 2, 10, 27, 28]. Этому планируется посвятить отдельную публикацию. Для более основательного знакомства с приложениями, в которых возникают задачи ЭЛП можно рекомендовать [3 – 9]. Особенно интересно сопоставить предложенные в данной статье методы с методом балансировки [1] (говорят также методом Шелейховского, Брэгмана–Шелейховского [9], Синхорна [28]) применительно к задаче расчета матрицы корреспонденций по энтропийной модели [3, 9, 17]:

$$f(x) = \sum_{i,j=1}^{n'} x_{ij} \ln x_{ij} + \alpha \sum_{i,j=1}^{n'} c_{ij} x_{ij} \to \min_{\substack{\sum_{j=1}^{n'} x_{ij} = L_i, \sum_{i=1}^{n'} x_{ij} = W_j \\ i,j=1,\dots,n; \, x \in S_{n'^2}(1)}}.$$

Если $\alpha \to \infty$ и $f(x) := f(x)/\alpha$, то рассматриваемая задача ЭЛП переходит в классическую транспортную задачу ЛП, которая в общем случае решается за $\mathrm{O}(n'^3 \ln n')$ арифметических операций, причем данная оценка не улучшаема [29]. Задача же ЭЛП может быть решена намного эффективнее, в частности, методом балансировки.

Для этой задачи, которую несложно привести к виду (1), двойственная задача будет иметь вид

$$\varphi(\lambda,\mu) = \langle \lambda, L \rangle + \langle \mu, W \rangle - \ln\left( \sum_{i,j=1}^{n'} \exp(-\alpha c_{ij} + \lambda_i + \mu_j) \right) \to \max_{\lambda, \mu \in \mathbb{R}^{n'}}, \qquad (8)$$

где $\sum_{i=1}^{n'} L_i = \sum_{j=1}^{n'} W_j = 1$. Метод балансировки примет вид ($[\lambda]_0 = [\mu]_0 = 0$):

$$[\lambda_i]_{k+1} = -\ln\left( \frac{1}{L_i} \sum_{j=1}^{n'} \exp(-\alpha c_{ij} + [\mu_j]_k) \right),$$

$$[\mu_j]_{k+1} = -\ln\left( \frac{1}{W_j} \sum_{i=1}^{n'} \exp(-\alpha c_{ij} + [\lambda_i]_k) \right) \text{ или } [\mu_j]_{k+1} = -\ln\left( \frac{1}{W_j} \sum_{i=1}^{n'} \exp(-\alpha c_{ij} + [\lambda_i]_{k+1}) \right).$$

Эти формулы можно получить, если заметить, что задача оптимизации (8) может быть явно решена по $\lambda$ при фиксированном $\mu$, и наоборот. Отметим при этом, что если $(\lambda, \mu)$ – решение задачи (8), то $(\lambda + c_1 e_{n'}, \mu + c_2 e_{n'})$, где $e_{n'}$ – вектор из $n'$ единиц, $c_1$, $c_2$ – произвольные числа, также будет решением задачи (8). Собственно, выписанные формулы – есть не что иное как, метод простой итерации для явно выписываемых условий экстремума (принципа Ферма) для задачи (8):

$$\lambda = \Lambda(\mu), \ \mu = \mathrm{M}(\lambda).$$

Оператор $(\lambda, \mu) \to (\Lambda(\mu), \mathrm{M}(\lambda))$ является сжимающим в метрике Биркгофа–Гильберта [30]. Причем можно оценить коэффициент сжатия, что в итоге приводит к следующим оценкам скорости сходимости (по числу итераций) метода балансировки

$$N = \mathrm{O}\left( \sqrt{\theta} \ln\left( \frac{\varepsilon_0}{\varepsilon} \right) \right),$$

где точность $\varepsilon$ – расстояние по метрике Биркгофа–Гильберта от того, что выдает метод до решения (неподвижной точки), $\varepsilon_0$ – расстояние от точки старта до неподвижной точки, а

$$\theta = \max_{i,j,p,q = 1,\ldots,n'} \exp\left( \alpha \cdot \left( c_{iq} + c_{pj} - (c_{ij} + c_{pq}) \right) \right).$$

Заметим, что параметр $\theta$ может быть содержательно проинтерпретирован исходя из эволюционного вывода модели расчета матрицы корреспонденций [9, 17].

Был проведен ряд численных экспериментов для решения задачи расчета матрицы корреспонденций (рис. 1 – 4). Эксперименты проводились на ЭВМ с процессором Intel Core i5, 2.5 ГГц и оперативной памятью 2 Гб в среде Matlab 2012® (8.0) под управлением операционной системой Microsoft Windows 7 (64 разрядная). Значительная часть параметров заполнялась автоматически: компоненты матрицы $\left\|c_{ij}\right\|_{i,j=1,1}^{n',n'}$ – случайные числа из интервала $(0,1)$; аналогично выбирались (а потом шкалировались) $L_i$, $W_j$, $\alpha = 100$ (на рис. 1, 2, 3). Требования к точности задавались относительной ошибкой в 1% (на рис. 1, 2, 4):

$$\varepsilon_f = 0.01 f(x_0), \; \varepsilon = 0.01 \|Ax_0 - b\|_2,$$

где $x_0 = x(\lambda = 0, \mu = 0)$ – точка старта.

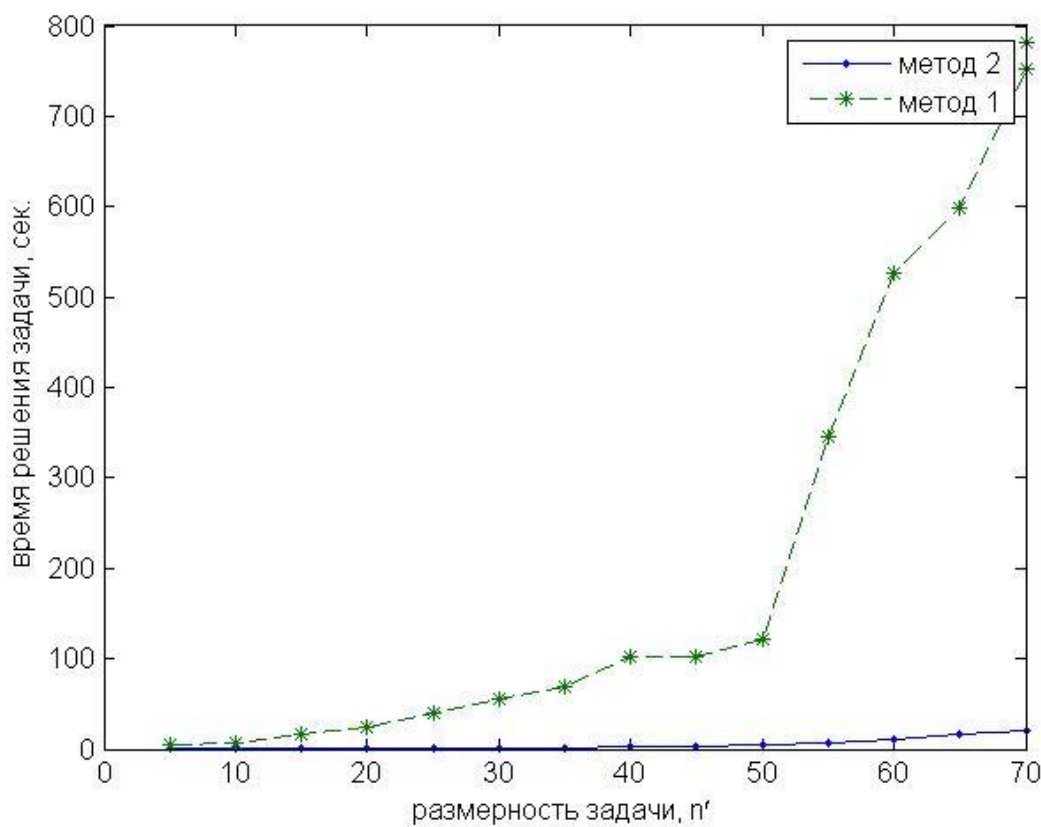

**Рис. 1. Сравнение зависимостей времен работ методов из теорем 1 и 2 от $n'$ (точность решения 1%)**

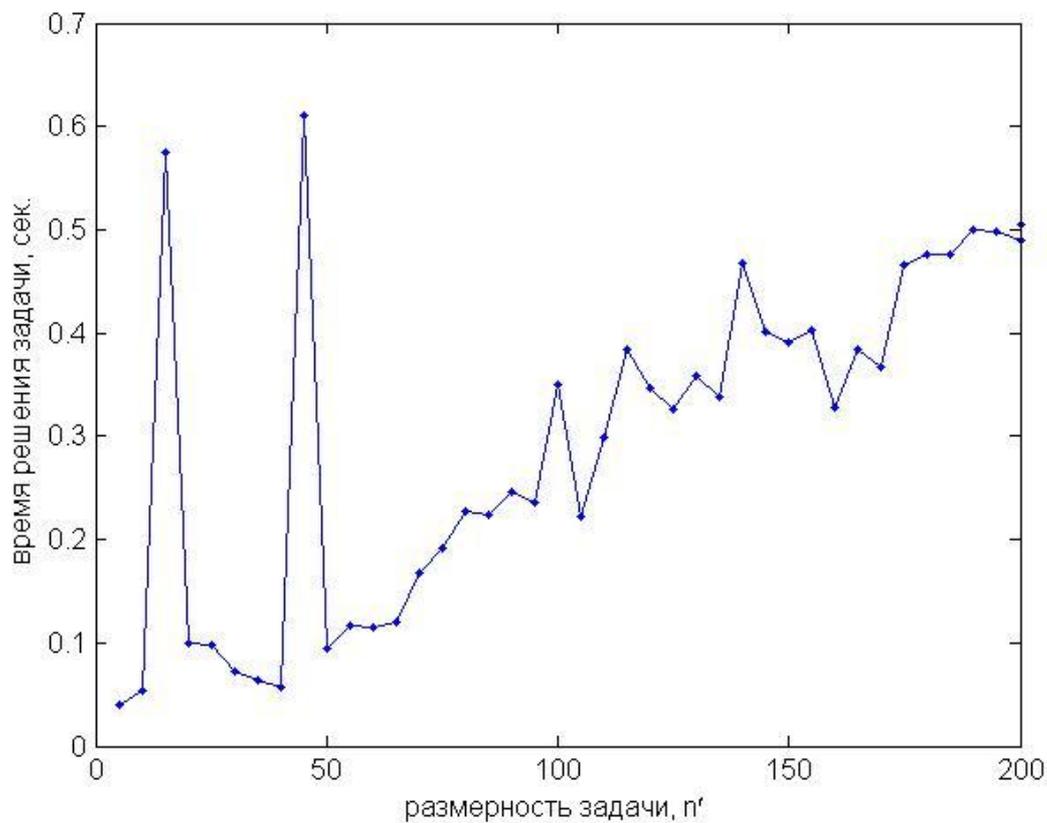

**Рис. 2. Зависимость времени работы метода балансировки от** $n'$ **(точность решения 1%)**

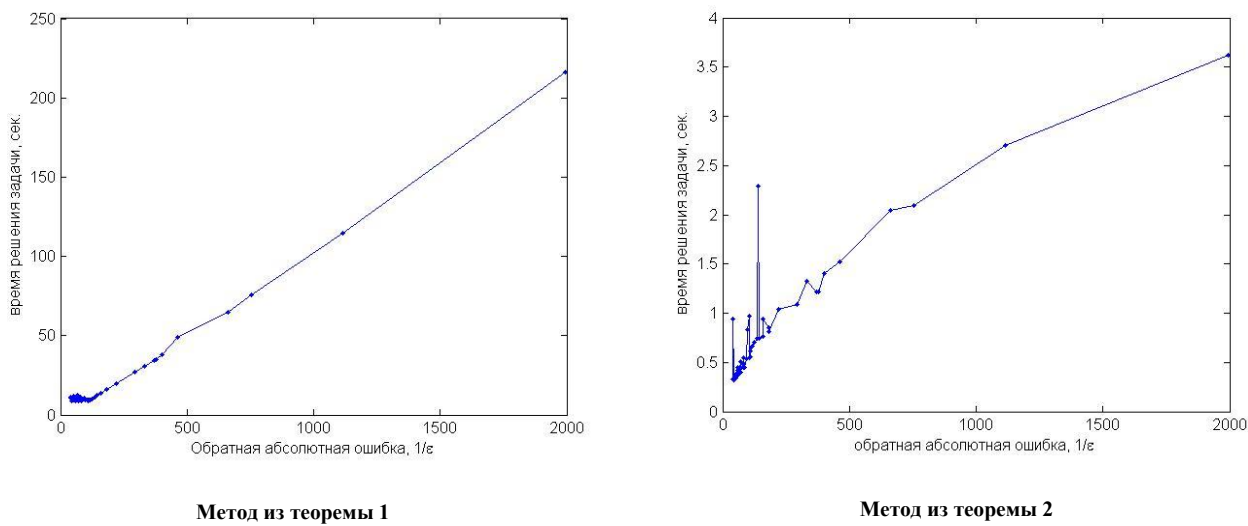

**Метод из теоремы 1**          **Метод из теоремы 2**

**Рис. 3 Зависимости времен работ методов из теорем 1 и 2 от точности** $\varepsilon$ ( $n' = 25$ )

График зависимости времени работы метода балансировки от точности (на той же задаче, что и на рис. 3) имеет вид горизонтальной прямой на уровне одной секунд, т.е. зависимость времени работы метода балансировки от точности крайне слабая (логарифмическая).

Количество итераций, которые делают методы при $\alpha = 100$, $n' = 30$ на относительной точности 1%, следующие: балансировка (58), метод из теоремы 2 (10 346), метод из теоремы 1 (221 481). Видно, что метод 1 делает много "лишних" итераций. Сложность одной итерации у методов 1 и 2 одинакова, но в методе 1 предписано сделать большее число итераций согласно полученной в теореме 1 оценке. Если бы в качестве критерия остановки в методе 1 мы взяли бы проверку условий леммы 1, то это немного увеличило бы стоимость одной итерации, но при этом число итераций стало бы приблизительно таким же, как и в методе 2. Численные эксперименты показывают, что в таком случае методы 1 и 2 требуют количества итераций, отличающиеся не более чем в 3 – 4 раза.

Интересно отметить, что время работы метода балансировки растет с ростом $\alpha$, в то время как время работы методов 1 и 2 с ростом $\alpha$ уменьшается (см. рис. 4) в наблюдаемом диапазоне (далее, при больших значениях $\alpha$, наблюдалась стабилизация).

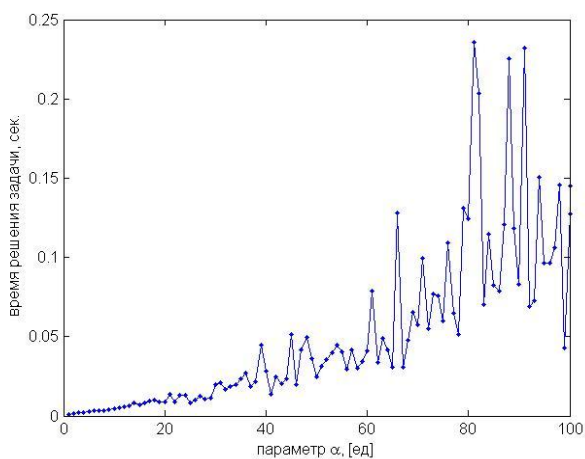
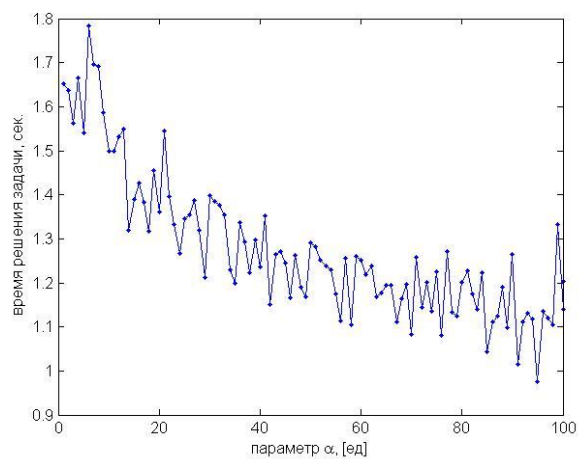

**Метод балансировки**            **Метод из теоремы 2**

**Рис. 4 Зависимости времен работ метода балансировки и метода из теорем 2 от $\alpha$ ( $n' = 30$ )**

Проведенные численные эксперименты показали, что предложенные в статье методы уступают по эффективности методу балансировки. Это не удивительно, поскольку балансировка считается наилучшим методам для специального класса задач ЭЛП (задач расчета матриц корреспонденций по энтропийной модели). Предложенные же в статье методы применимы для более широкого класса задач ЭЛП, в том числе с ограничениями в виде линейных неравенств (и даже более общих конусных неравенств вида $Ax - b \in K$, с двойственным конусом $K^*$ просто структуры). Также результаты статьи не сложно перенести на случай, когда вместо единичного симплекса рассматривается множество

$$\left\{ x \in \mathbb{R}^n : x_i \geq 0, \ i = 1,...,n, \ \sum_{i=1}^{n} x_i \leq 1 \right\}.$$

## Список литературы

1. *Fang S.-C., Rajasekera J.R., Tsao H.-S.J.* Entropy optimization and mathematical programming. Kluwer's International Series, 1997.

2. *Попков Ю.С.* Теория макросистем: Равновесные модели. М.: УРСС, 2013.

3. *Вильсон А.Дж.* Энтропийные методы моделирования сложных систем. М.: Наука, 1978.

4. International workshops on Bayesian inference and maximum entropy methods in science and engineering. AIP Conf. Proceedings (holds every year from 1980).

5. *Kapur J.N.* Maximum – entropy models in science and engineering. John Wiley & Sons, Inc., 1989.

6. *Golan A., Judge G., Miller D.* Maximum entropy econometrics: Robust estimation with limited data. Chichester, Wiley, 1996.

7. *Jaynes E.T.* Probability theory. The logic of science. Cambridge: Cambridge University Press, 2003.

8. *Гасников А.В.* Марковские модели макросистем // e-print, 2014. arXiv:1412.2720

9. *Гасников А.В., Кленов С.Л., Нурминский Е.А., Холодов Я.А., Шамрай Н.Б.* Введение в математическое моделирование транспортных потоков. Под ред. А.В. Гасникова с приложениями М.Л. Бланка, К.В. Воронцова и Ю.В. Чеховича, Е.В. Гасниковой,

А.А. Замятина и В.А. Малышева, А.В. Колесникова, Ю.Е. Нестерова и С.В. Шпирко, А.М. Райгородского, с предисловием руководителя департамента транспорта г. Москвы М.С. Ликсутова. М.: МЦНМО, 2013. 427 стр., 2-е изд.

http://www.mou.mipt.ru/gasnikov1129.pdf

10. *Гасникова Е.В.* Двойственные мультипликативные алгоритмы для задач энтропийно-линейного программирования // ЖВМ и МФ. 2009. Т.49 № 3. С. 453–464.

11. *Немировский А.С., Юдин Д.Б.* Сложность задач и эффективность методов оптимизации. М.: Наука, 1979. http://www2.isye.gatech.edu/~nemirovs/Lect_EMCO.pdf

12. *Кельберт М.Я. , Сухов Ю.М.* Марковские цепи как отправная точка теории случайных процессов и их приложений. Вероятность и статистика в примерах и задачах. Т. 2. М.: МЦНМО, 2010.

13. *Кац М.* Вероятность и смежные вопросы в физике. М.: Мир, 1965.

14. *Малышев В.А., Пирогов С.А.* Обратимость и необратимость в стохастической химической кинетике // Успехи мат. наук. 2008. Т. 63. вып. 1(379). С. 4–36.

15. *Гасников А.В., Гасникова Е.В.* Об энтропийно-подобных функционалах, возникающих в стохастической химической кинетике при концентрации инвариантной меры и в качестве функций Ляпунова динамики квазисредних // Матем. Заметки. 2013. Т. 94:6. С. 819–827. arXiv:1410.3126

16. *Веденяпин В.В. , Аджиев С.З.* Энтропия по Больцману и Пуанкаре // Успехи мат. наук. 2014. Т. 69:6(420). С. 45–80.

17. *Гасников А.В., Гасникова Е.В., Мендель М.А., Чепурченко К.В.* Эволюционные выводы энтропийной модели расчета матрицы корреспонденций // Мат. мод. 2016. Т. 28. (принята к печати) arXiv:1508.01077

18. *Nesterov Y.* Smooth minimization of non-smooth function // Math. Program. Ser. A. 2005. V. 103. № 1. P. 127–152.

19. *Нестеров Ю.Е.* Введение в выпуклую оптимизацию. М.: МЦНМО, 2010.

20. *Гасникова Е.В.* Моделирование динамики макросистем на основе концепции равновесия. Дисс. на соиск. степени к.ф.-м.н., спец. 05.13.18, 17 декабря 2012 г. 90 С.

21. *Devolder O., Glineur F., Nesterov Y.* Double smoothing technique for large-scale linearly constrained convex optimization // SIAM Journal of Optimization. 2012. V. 22. № 2. P. 702–727.

22. *Гасников А.В., Дорн Ю.В., Нестеров Ю.Е, Шпирко С.В.* О трехстадийной версии модели стационарной динамики транспортных потоков // Математическое моделирование. 2014. Т. 26:6. С. 34–70. arXiv:1405.7630

23. *Гасников А.В., Двуреченский П.Е., Нестеров Ю.Е.* Стохастические градиентные методы с неточным оракулом // Труды МФТИ. 2016. Т. 8. arxiv:1411.4218

24. *Juditsky A., Nemirovski A.* First order methods for nonsmooth convex large-scale optimization, II: Utilizing problem's structure. In: Optimization for Machine Learning. Eds. S. Sra, S. Nowozin, S. Wright. MIT Press, 2012.

    http://www2.isye.gatech.edu/~nemirovs/MLOptChapterII.pdf

25. *Bubeck S.* Convex optimization: algorithms and complexity // In Foundations and Trends in Machine Learning. 2015. V. 8. no. 3-4. P. 231–357. arXiv:1405.4980

26. *Назин А.В., Тремба А.А.* Игровой алгоритм зеркального спуска в задаче робастного PageRank // Автоматика и телемеханика. 2016. (принята к печати)

27. *Boyd S. et al.* Open Software http://www.cvxpy.org/en/latest/examples/index.html ; http://nbviewer.ipython.org/github/cvxgrp/cvxpy/blob/master/examples/notebooks/WWW/max_entropy.ipynb

28. *Cuturi M.* Personal page http://www.iip.ist.i.kyoto-u.ac.jp/member/cuturi/

29. *Pele O., Werman M.* Fast and robust earth mover's distances. In ICCV'09, 2009.

    http://www.cs.huji.ac.il/~werman/Papers/ICCV2009.pdf

30. *Franklin J., Lorenz J.* On the scaling of multidimensional matrices // Linear Algebra and its applications. 1989. V. 114. P. 717–735.